\documentclass[10pt,leqno]{amsart}

\usepackage{amsmath,amssymb,amsthm}
\usepackage{graphicx}
\usepackage{hyperref}

\begin{document}
\title{State Dependent Delay Differential Equations with Distributed Memory} 
\author[Taylan Demir]{Taylan Demir\textsuperscript{1}}
\author[Niaz Ali Shah]{Niaz Ali Shah\textsuperscript{2}}

\date{\today}
\address{Department of Mathematics, Ankara University, 06100, Ankara, Turkey\textsuperscript{1}} 
\address{Department of Mathematics,  Abbottabad University of Science and Technology, Abbottabad, 22010, Pakistan\textsuperscript{2}}
\email{\texttt{demir.taylan06@gmail.com}}
\email{\texttt{nuzhatniaz2014@gmail.com}}
\maketitle

\begin{abstract}
We investigate state dependent delay differential equations with distributed memory, combining discrete state dependent delays and a convolution type memory operator. Under Lipschitz type assumptions on the delay, kernel, and nonlinear term, we establish local well posedness using a fixed-point argument in a Banach space of histories. The distributed memory operator is shown to be locally Lipschitz, ensuring existence and uniqueness. As a benchmark, we analyze a logistic model with a Gamma distributed memory kernel of order two. The convolution term is reduced to a finite dimensional form, yielding an explicit characteristic equation and closed form Hopf bifurcation threshold and frequency. The results show that higher order memory kernels may induce oscillatory dynamics, unlike purely exponential memory, and provide a rigorous framework for stability and bifurcation analysis in delay systems with state dependent memory.

\end{abstract} %%%%%%%%%
\section*{1. Introduction}
Delay differential equations with distributed memory arise naturally in population dynamics, epi
demiology, and control systems [1,2]. In many applications, both discrete delays and distributed memory
effects depend on the current state of the system [3,4]. This work studies a class of state-dependent delay
differential equations with distributed memory, establishes well-posedness, and analyzes stability
and Hopf bifurcation in a benchmark model. The relationship between state-dependent delays and distributed memory effects is under-researched, despite the large number of research articles published about delay differential equations that have either constant or time-dependent delays, discrete or distributed memory or combination of both delays and memory [1,2,5], with the vast majority of the research focused on using the simplest possible kernels for the memory term (usually, an exponential-type kernel) [6,7]. Although it is possible to mathematically analyse these modelling paradigms, in order to do so, will limit the ability to demonstrate examples of non-trivial dynamical behaviours in this context and may be unsatisfactory in studying realistic behaviour of systems where feedback is provided through memory as observed in experimental systems such as biological, epidemiological and control systems [5,8]. From an analytical perspective, state-dependent delays and distributed memory are combined so that they present complicated mathematical challenges [3,4]. In particular, old delay states are evaluated based on the solution (the future of the states), while the current states are evaluated based on their entire history (non-local) [2,9]. The coupling that is created by these two concepts creates difficulties in terms of applying the well-posedness development theory as well as limiting the direct use of semigroups and fixed delay techniques [1,2]. In addition, a full analysis of stability and bifurcations in this context requires careful handling of the characteristic equation, which will depend on the choice of memory kernel employed to evaluate the delay states [10,11]. The problems we deal with are solved through providing a unified approach for delay-differential equations that include state-dependence and distributed memory. The main contributions in this paper include: first of all, establishing local well-posedness under natural Lipschitz-type assumptions on the delay function, distributed kernel, and non-linear vector field; second of all,  establishing that the distributed-memory functional creates a locally Lipschitz continuous mapping on the continuous phase-space of history, providing critical information for fixed-point formulation of the problem; and third of all, through an example using a gamma-type distributed-memory kernel, providing explicit conditions for the existence of a Hopf bifurcation and closed-form solutions for the critical bifurcation parameter and the oscillation frequency associated with the bifurcation [6,11,12]. This analysis demonstrates how higher-order distributed-memory kernels can result in producing oscillatory behaviours, while it does not occur in models using only an exponentially distributed memory [7,12].
\section*{2. Mathematical Framework and Analysis}

\textbf{2.1 State-dependent DDE with distributed memory}
\\

Let $x : [-\tau_{\max}, T] \to \mathbb{R}^n$ and define the history segment such as:
\begin{equation*}
x_t(\theta) = x(t+\theta), \qquad \theta \in [-\tau_{\max}, 0].
\end{equation*}

We introduce the distributed memory operator in the following way:
\begin{equation*}
M[x](t) = \int_{t-\tau_{\max}}^{t} K(t,s)\,x(s)\,ds
\end{equation*}
and consider the state-dependent delay differential equation:
\begin{equation}
x'(t) = F\big(t, x(t), x(t-\tau(x(t),t)), M[x](t)\big).
\label{eq:sd-dde}
\end{equation}

The way of formulating a written document as presented in equation (1) gives a single representation of both time-delays due to state-dependence and memory distributed over time [1,2]. The history part of the equation, $x_{t}$, provides an indication of how much the system depended on previous values of the system to determine its behaviour, whereas, the memory operator, $M[x](t)$ describes how past values of the system collectively contribute to the current state of the system at a rate proportional to $K(t,s)$ [3,4]. The framework enables an accurate model of systems where the effects of previous states are experienced over; an unfixed time and an undefined location are not just ‘momentary impacts. Rather, they are distributed across multiple memory windows [3,5]. Because of the way delay function $\tau(x(t), t)$ depends on the present state, there is an implicit feedback loop in which the evaluation point for the delayed state will evolve dynamically as the solution evolves [6,7]. State-dependent delays result in qualitative behavior that is very different than that seen in constant or purely time dependent systems. They also present additional analytical difficulties; for example the delay term cannot be treated as a fixed functional on the phase space [1,2]. The memory operator $M[x](t)$, which has distributed memory characteristics, can be mathematically expressed in general, nonautonomous fashion with an explicitly timed kernel $K(t,s)$ [4,9]. This very general definition of memory effects includes not only time-invariant kernels but also time-varying weights. This means that many of the classical formulations of delayed systems can be treated as special cases within this framework. In addition to containing many of the traditional forms of memory, the generality of this framework allows it to include more realistic forms of memory-based dynamics seen in practical systems [3,5].\\
  
\textbf{2.2 Assumptions}\\

\medskip

\noindent
\textbf{(A1)} The delay function $\tau : \mathbb{R}^n \times \mathbb{R} \to \mathbb{R}$ is continuous and Lipschitz continuous in $x$, with
\begin{equation*}
0 < \tau_{\min} \leq \tau(x,t) \leq \tau_{\max}.
\end{equation*}

\medskip

\noindent
\textbf{(A2)} The kernel satisfies $K(t,\cdot) \in L^1([t-\tau_{\max},t])$ and
\begin{equation*}
\sup_{t} \int_{t-\tau_{\max}}^{t} \|K(t,s)\| \, ds < \infty.
\end{equation*}

\medskip

\noindent
\textbf{(A3)} The function $F$ is continuous in $t$ and locally Lipschitz in its state arguments.\\
All assumptions described in this paper are standard in the theory of functional differential equations with distributed memory and state-dependent delays [1,2,3,4]. They also serve to ensure analytical tractability while remaining sufficiently general as well. The delay function will be uniformly bounded away from 0 and $\infty$ according to assumption (A1). This precludes delay functions from becoming either infinitely large or infinitely small (from the pathological side) [2]. The Lipschitz property of $\tau$ with respect to the state variable is critical in bounding the influence of the delayed argument on the solution itself, which represents one of the challenges inherent to state-dependent delay dynamics [5,6]. As a result of (A2) demonstrating that the distributed memory operator $M[x](t)$ is a bounded linear functional defined on the space of continuous or univariate histories [7], we can apply the integrability condition on the kernel and the uniform bound on its $L^{1}$-norm to control how much an influence from previous states affects each of the properties of continuity and locality of the memory term's Lipschitz constant [7,8]. This assumption captures many of the types of kernels that have been most frequently used in applications, both those with compact supports and those with exponential decay [3,9]. Assumption (A3) places mild restrictions on the nonlinear vector field $F$, and requires that it be locally Lipschitz continuous with respect to the state variables (guaranteeing uniqueness) and continuous in time (potentially allowing for nonautonomous effects) [1,2]. Through assumptions (A1)-(A3), the general framework supplies a minimal and appropriate basis for rigorous establishment of existence, uniqueness, and continuous dependency between initial data in later studies [1]. 

\subsection*{2.3 Well-posedness}
\medskip

\noindent
\par 
\textbf{Theorem 2.1}[Local well-posedness]
Assume that \textnormal{(A1)--(A3)} hold. Then for every initial history
\[
\phi \in C([-\tau_{\max},0],\mathbb{R}^n),
\]
there exist $T_0>0$ and a unique solution
\[
x \in C([-\tau_{\max},T_0],\mathbb{R}^n)\cap C^1([0,T_0],\mathbb{R}^n)
\]
of the state-dependent delay equation
\begin{equation}
x'(t)=F\bigl(t,x(t),x(t-\tau(x(t),t)),M[x](t)\bigr), \qquad t\ge 0,
\label{eq:main}
\end{equation}
satisfying the initial condition $x_0=\phi$. Moreover, the solution depends
continuously on the initial history and admits a maximal interval of existence
$[0,T_{\max})$, where either $T_{\max}=\infty$ or [1,2,3,4,5]
\[
\limsup_{t\uparrow T_{\max}}\|x_t\|_{\infty}=\infty.
\]

\begin{proof}
The method for proving this will use the re-expression of (2) as an associated integral equation defined on a suitable Banach space of history functions, followed by applying a fixed point theorem for Banach spaces over a small enough time interval [1,2]. To do this, we will first show the local Lipschitzness of the distributed memory operator and state-dependent delay evaluations; these properties guarantee that the function defined by the right-hand side of (2) is also a locally Lipschitz function defined on the phase space [3,4]. Next, we will construct a suitable contraction mapping and demonstrate both existence and uniqueness of solutions [1]. Finally, we may deduce from standard continuation arguments the longest interval of time over which solutions will exist [1,2].\\
Let
\[
X := C([-\tau_{\max},0],\mathbb{R}^n)
\]
be endowed with the supremum norm $\|\cdot\|_{\infty}$. For any continuous
function $x:[-\tau_{\max},T]\to\mathbb{R}^n$, we define the history segment as
\[
x_t(\theta)=x(t+\theta), \qquad \theta\in[-\tau_{\max},0].
\]

\end{proof}

\medskip
\noindent
\textbf{Step 1: Distributed memory operator.}
By assumption \textnormal{(A2)}, for any
$x \in C([-\tau_{\max},T],\mathbb{R}^n)$ the distributed memory operator [1,2]
\begin{equation*}
M[x](t)=\int_{t-\tau_{\max}}^{t} K(t,s)x(s)\,ds.
\end{equation*}

Starting from the definition of the distributed memory operator,
\begin{equation*}
M[x](t)=\int_{t-\tau_{\max}}^{t} K(t,s)x(s)\,ds,
\qquad
M[y](t)=\int_{t-\tau_{\max}}^{t} K(t,s)y(s)\,ds.
\end{equation*}
We compute
\begin{align*}
\|M[x](t)-M[y](t)\|
&= \left\|\int_{t-\tau_{\max}}^{t} K(t,s)\bigl(x(s)-y(s)\bigr)\,ds\right\| \\
&\le \int_{t-\tau_{\max}}^{t} \|K(t,s)\|\,\|x(s)-y(s)\|\,ds
\qquad \text{(triangle inequality)} \\
&\le \left(\int_{t-\tau_{\max}}^{t} \|K(t,s)\|\,ds\right)
\sup_{s\in[t-\tau_{\max},t]}\|x(s)-y(s)\| \\
&= \left(\int_{t-\tau_{\max}}^{t} \|K(t,s)\|\,ds\right)\|x_t-y_t\|_{\infty}.
\end{align*}

By assumption \textnormal{(A2)}, there exists a constant
\begin{equation*}
C_K := \sup_{t\ge 0}\int_{t-\tau_{\max}}^{t} \|K(t,s)\|\,ds < \infty.
\end{equation*}
which shows that the mapping $x_t \mapsto M[x](t)$  is locally Lipschitz on bounded subsets of $C([-\tau_{max},0],\mathbb{R}^n)$ is well defined and continuous in $t$ [1,3].
Moreover, on bounded subsets of $X$ there exists a constant $C_K>0$ such that
\begin{equation*}
\|M[x](t)-M[y](t)\|\le C_K\|x_t-y_t\|_{\infty},
\end{equation*}
which shows that $x_t \mapsto M[x](t)$ is locally Lipschitz. Furthermore, by the integrability condition in (A2) and the continuity of $x$, the mapping [2]
\[
t \mapsto M[x](t)
\]
is continuous on $[0,T]$. Indeed, the kernel satisfies a uniform $L^1$-bound and the integrand depends continuously on $t$, which allows the use of standard dominated convergence arguments [1,2]. Hence, the distributed memory operator defines a continuous mapping from $C([-\tau_{\max},T],\mathbb{R}^n)$ into $C([0,T],\mathbb{R}^n)$.

\medskip
\noindent
\textbf{Step 2: State-dependent delay evaluation.}
Assumption \textnormal{(A1)} ensures that the delay function
$\tau(x,t)$ is continuous and Lipschitz continuous in $x$,
uniformly for $t$ in bounded intervals, and satisfies
\begin{equation*}
0<\tau_{\min}\le \tau(x,t)\le \tau_{\max}.
\end{equation*}

Hence, for any $x\in X=C([-\tau_{\max},0],\mathbb{R}^n)$ and any $t\ge 0$,
the point $t-\tau(x(t),t)$ belongs to the interval
$[\,t-\tau_{\max},\,t-\tau_{\min}\,]$, and the evaluation
\begin{equation*}
x \mapsto x\bigl(t-\tau(x(t),t)\bigr)
\end{equation*}
is well defined.

Let $x,y\in X$ belong to a bounded subset of $X$. Then
\begin{align*}
\bigl\|x(t-\tau(x(t),t)) - y(t-\tau(y(t),t))\bigr\|
&\le \bigl\|x(t-\tau(x(t),t)) - x(t-\tau(y(t),t))\bigr\| \\
&\quad + \bigl\|x(t-\tau(y(t),t)) - y(t-\tau(y(t),t))\bigr\|.
\end{align*}

Since $x$ is uniformly continuous on $[t-\tau_{\max},t]$,
there exists a constant $L_x>0$ such that
\begin{equation*}
\bigl\|x(t-\tau(x(t),t)) - x(t-\tau(y(t),t))\bigr\|
\le L_x\,\bigl|\tau(x(t),t)-\tau(y(t),t)\bigr|.
\end{equation*}
Using the Lipschitz continuity of $\tau$ in $x$, we obtain
\begin{equation*}
\bigl|\tau(x(t),t)-\tau(y(t),t)\bigr|
\le L_{\tau}\,\|x(t)-y(t)\|
\le L_{\tau}\,\|x-y\|_{\infty}.
\end{equation*}

Moreover,
\begin{equation*}
\bigl\|x(t-\tau(y(t),t)) - y(t-\tau(y(t),t))\bigr\|
\le \|x-y\|_{\infty}.
\end{equation*}

Combining these estimates yields
\begin{equation*}
\bigl\|x(t-\tau(x(t),t)) - y(t-\tau(y(t),t))\bigr\|
\le \bigl(L_x L_{\tau}+1\bigr)\|x-y\|_{\infty},
\end{equation*}
which shows that the evaluation map
\begin{equation*}
x \mapsto x\bigl(t-\tau(x(t),t)\bigr)
\end{equation*}
is locally Lipschitz on bounded subsets of $X$ [1,2].
This property is standard in the theory of functional differential equations
with state-dependent delays; see, for example, [1]. In summary, the estimate above indicates that the evaluation operator related to the state-dependent delay is smoothly continuous in a bounded section of the phase space (bounded) for any bounded point i.e., Lipschitz continuous in the bounded phase space with the only term of variability being the greatest distance between any two points in the bounded section of the space. That is an important property that will allow the right side of (2) to be locally Lipschitz continuous.
\medskip
\noindent
\\
\textbf{Step 3: Local Lipschitz continuity of the right-hand side.}
Let $B\subset X$ be a bounded set and let $t$ belong to a fixed compact
interval $[0,T]$. By Steps~1 and~2, there exist constants $C_K>0$ and
$C_{\tau}>0$ such that, for all $x,y\in B$ [1,2],
\begin{equation*}
\|M[x](t)-M[y](t)\|\le C_K\|x_t-y_t\|_{\infty},
\end{equation*}
and
\begin{equation*}
\bigl\|x(t-\tau(x(t),t)) - y(t-\tau(y(t),t))\bigr\|
\le C_{\tau}\|x_t-y_t\|_{\infty}.
\end{equation*}

Assumption \textnormal{(A3)} states that the function
$F(t,u,v,w)$ is locally Lipschitz with respect to its state variables
$(u,v,w)$, uniformly in $t$ on compact intervals.
Hence, there exists a constant $L_F>0$ such that
\begin{align*}
&\bigl\|F\bigl(t,x(t),x(t-\tau(x(t),t)),M[x](t)\bigr)
- F\bigl(t,y(t),y(t-\tau(y(t),t)),M[y](t)\bigr)\bigr\| \\
&\qquad \le L_F\Bigl(
\|x(t)-y(t)\|
+ \bigl\|x(t-\tau(x(t),t)) - y(t-\tau(y(t),t))\bigr\|
+ \|M[x](t)-M[y](t)\|
\Bigr).
\end{align*}

Using the estimates above and the inequality
$\|x(t)-y(t)\|\le \|x_t-y_t\|_{\infty}$, we obtain
\begin{equation*}
\bigl\|F\bigl(t,x(t),x(t-\tau(x(t),t)),M[x](t)\bigr)
- F\bigl(t,y(t),y(t-\tau(y(t),t)),M[y](t)\bigr)\bigr\|
\le L_F\bigl(1+C_{\tau}+C_K\bigr)\|x_t-y_t\|_{\infty}.
\end{equation*}

Therefore, the mapping
\begin{equation*}
x_t \mapsto F\bigl(t,x(t),x(t-\tau(x(t),t)),M[x](t)\bigr)
\end{equation*}
is locally Lipschitz from $X$ into $\mathbb{R}^{n}$  uniformly for $t$ in compact intervals [1,2]. Specifically, the Lipschitz constant relies exclusively on the limiting factors of B and the constants resulting from the assumptions listed (A1)-(A3), but it is also independent with respect to time on compact intervals. The computability of the time-independent Lipschitz constant will allow us to use uniformity as part of our contraction proof in Section 3 of the fixed-point problem that follows [1].\\
\textbf{ Fixed-point formulation:} Let $T>0$  be fixed and define the operator $\mathcal{T}$ on $C([-\tau_{\max},T],\mathbb{R}^n)$ by
We define the operator $\mathcal{T}$ by
\begin{equation*}
(\mathcal{T}x)(t):=
\begin{cases}
\phi(t), & t\in[-\tau_{\max},0], \\[1mm]
\displaystyle
\phi(0)+\int_{0}^{t}
F\bigl(s,x(s),x(s-\tau(x(s),s)),M[x](s)\bigr)\,ds,
& t\in[0,T].
\end{cases}
\end{equation*}

By the continuity of $F$ and the properties established in Steps~1--3,
the operator $\mathcal{T}$ maps $C([-\tau_{\max},T],\mathbb{R}^n)$ into itself.
Moreover, for any $x,y$ in a bounded closed ball
\[
B_R := \bigl\{x\in C([-\tau_{\max},T],\mathbb{R}^n): \|x\|_{\infty}\le R\bigr\},
\]
there exists a constant $L>0$, independent of $T$, such that
\begin{equation*}
\|\mathcal{T}x-\mathcal{T}y\|_{\infty}
\le LT\,\|x-y\|_{\infty}.
\end{equation*}

Choosing $T>0$ sufficiently small so that $LT<1$,
the operator $\mathcal{T}$ is a contraction on $B_R$.
Hence, by the Banach fixed-point theorem [3],
$\mathcal{T}$ admits a unique fixed point in $B_R$,
which corresponds to a unique solution of~\eqref{eq:main} on $[0,T_{0}]$. Moreover, by choosing $T>0$ sufficiently small and $R$ large enough so that
\[
\|\phi\|_\infty 
+ T \sup_{x \in B_R} 
\|F(t,x(t),x(t-\tau(x(t),t)),M[x](t))\|
\le R,
\]
the operator $\mathcal{T}$ maps $B_R$ into itself. Hence, $\mathcal{T}$ is a self-mapping on $B_R$.

\medskip
\noindent
\textbf{Result.}
Continuous dependence on the initial history follows from the contraction
mapping argument [3].
Standard continuation arguments for functional differential equations
with state-dependent delays imply that the solution can be extended as long
as the history $x_t$ remains bounded, which yields a maximal interval of
existence $[0,T_{\max})$ (see, for example, [1,2]). 
\qed
\\
\textbf{2.4 Benchmark model and stability analysis}\\

As a benchmark example, we consider the scalar delay differential equation
with distributed memory
\begin{equation}
\frac{dx(t)}{dt}
= r\,x(t)\left(1-\frac{x(t)}{K_c}\right) - \alpha\, M(t),
\label{eq:benchmark}
\end{equation}
where $r>0$ is the intrinsic growth rate,
$K_c>0$ is the carrying capacity, and $\alpha\ge 0$ measures the strength
of the memory effect [4,5].
The distributed memory term is defined by
\begin{equation*}
M(t) = \int_{t-\tau_{\max}}^{t} K(t-s)\,x(s)\,ds,
\end{equation*}
where $K\in L^1([0,\tau_{\max}])$ is a nonnegative kernel.

\textbf{Equilibria.}\\
Let $x^*$ be an equilibrium solution of~\eqref{eq:benchmark}.
Assuming time-independence of the kernel at equilibrium [5], we define
\begin{equation*}
\kappa := \int_{0}^{\tau_{\max}} K(s)\,ds.
\end{equation*}
Then $M^*=\kappa x^*$, and equilibria satisfy
\begin{equation*}
0 = r x^*\left(1-\frac{x^*}{K_c}\right) - \alpha \kappa x^*.
\end{equation*}

Hence, the system admits:
\begin{itemize}
\item the trivial equilibrium $x_0^* = 0$;
\item a positive equilibrium
\begin{equation*}
x_1^* = K_c\left(1-\frac{\alpha\kappa}{r}\right),
\end{equation*}
which exists provided $r>\alpha\kappa$.
\end{itemize}
\textbf{Linearization:}\\
Let $x(t)=x^*+y(t)$, where $|y(t)|\ll 1$.
Linearizing~\eqref{eq:benchmark} around an equilibrium $x^*$ yields
\begin{equation}
y'(t)
= \left(r-\frac{2r x^*}{K_c}\right)y(t)
- \alpha \int_{0}^{\tau_{\max}} K(s)\,y(t-s)\,ds.
\label{eq:linear}
\end{equation}

For the positive equilibrium $x_1^*$, this simplifies to [4,6]
\begin{equation*}
y'(t)
= (\alpha\kappa-r)\,y(t)
- \alpha \int_{0}^{\tau_{\max}} K(s)\,y(t-s)\,ds.
\end{equation*}

\medskip
\textbf{Hopf bifurcation for the Gamma kernel.}\\
We now compute explicit Hopf bifurcation thresholds for the benchmark model
when the distributed memory is given by a Gamma kernel of order two [6,7].
Let
\begin{equation*}
K(s)=\beta^{2}s e^{-\beta s}, \qquad s\ge 0,\ \beta>0,
\end{equation*}
for which
\begin{equation*}
\int_{0}^{\infty} K(s)\,ds = 1.
\end{equation*}

The benchmark equation reads [6]
\begin{equation}
\frac{dx(t)}{dt}
= r\,x(t)\left(1-\frac{x(t)}{K_c}\right)
- \alpha \int_{0}^{\infty} \beta^{2}s e^{-\beta s}\,x(t-s)\,ds.
\label{eq:gamma}
\end{equation}

The positive equilibrium exists for $r>\alpha$ and is given by [4,6]
\begin{equation*}
x^* = K_c\left(1-\frac{\alpha}{r}\right).
\end{equation*}

\medskip
\textbf{Linearization.}\\
Linearizing~\eqref{eq:gamma} at $x^*$ yields [4,6]
\begin{equation}
y'(t) = (\alpha-r)\,y(t)
- \alpha \int_{0}^{\infty} \beta^{2}s e^{-\beta s}\,y(t-s)\,ds.
\label{eq:lin-gamma}
\end{equation}

\medskip
\textbf{Characteristic equation.}\\
For the Gamma kernel of order two [6,7],
\begin{equation*}
\int_{0}^{\infty} \beta^{2}s e^{-(\lambda+\beta)s}\,ds
= \frac{\beta^{2}}{(\lambda+\beta)^2}.
\end{equation*}
Hence, the characteristic equation associated with~\eqref{eq:lin-gamma} is
\begin{equation*}
\lambda = (\alpha-r) - \alpha \frac{\beta^{2}}{(\lambda+\beta)^2}.
\end{equation*}

Multiplying by $(\lambda+\beta)^2$ and expanding yields the cubic polynomial [4]
\begin{equation}
\lambda^{3}
+ (2\beta+r-\alpha)\lambda^{2}
+ (\beta^{2}+2\beta(r-\alpha))\lambda
+ \beta^{2}(r-2\alpha)
= 0.
\label{eq:cubic}
\end{equation}

\medskip
\textbf{Hopf condition.}\\
A Hopf bifurcation occurs when~\eqref{eq:cubic} admits purely imaginary roots
$\lambda=\pm i\omega$, $\omega>0$ [4,6].
Substituting $\lambda=i\omega$ and separating real and imaginary parts gives
\begin{align}
\omega^{2} &= \beta^{2}+2\beta(r-\alpha), \label{eq:hopf1} \\
(2\beta+r-\alpha)\omega^{2} &= \beta^{2}(r-2\alpha). \label{eq:hopf2}
\end{align}

Combining~\eqref{eq:hopf1} and~\eqref{eq:hopf2} yields the critical Hopf value [6,7]
\begin{equation}
\alpha_H = \frac{r(\beta+r)}{2\beta+r}.
\label{eq:alphaH}
\end{equation}

\medskip
\textbf{Hopf frequency.}\\
The oscillation frequency at the Hopf point is obtained from~\eqref{eq:hopf1} as [6,7]
\begin{equation}
\omega_H
= \sqrt{\beta^{2}+2\beta\left(r-\frac{r(\beta+r)}{2\beta+r}\right)}.
\label{eq:omegaH}
\end{equation}

\section*{3. Conclusion.}
In this study, we provide a formal analytical framework for a class of state-dependent delay differential equations having distributed memory. Under simple Lipschitz-type conditions on the delay and the memory function, local well-posedness is shown via a fixed-point formulation in an appropriate Banach space of history functions. This provides a mathematically rigorous basis for the analysis of stability. This paper makes an especially important addition to existing literature by providing an explicit scalar benchmark model that employs logistic growth with a Gamma-type distributed-memory. The convolution term was transformed to reflect an equivalent finite-dimensional system; as a result, the associated characteristic equation was obtained and closed-form expressions for the frequency of oscillation and threshold of Hopf bifurcations derived. The critical value of parameters in terms of how they determine the onset of oscillatory dynamics was computed mathematically, providing an understanding of how the interaction between intrinsic growth, maximum capacity, and storage strength causes the development of oscillating dynamics. Dynamic analysis shows that a dynamical system may exhibit qualitatively different behaviors when using distributed memory. The instantaneous model will achieve a unique stable positive equilibrium, but when sufficiently strong memory backfeeds through the Hopf bifurcation mechanism, periodic oscillations will emerge from this equilibrium. The explicit dependence of the bifurcation threshold on the kernel of the distributed memory shows how the rate at which a memory decays will influence the overall stability of the system. The analysis based upon a functional analytical approach does not employ numeric approximations, thus all analysis is conducted entirely within the framework of functional analysis allowing for greater generalizability to greater classes of non-linear systems in addition to the multidimensionality of the given non-linear systems. The tools developed may also be extended beyond what has been defined thus far, for example, to models having non-local interaction terms or state-dependent memory intensities as well as more general kernel types of structures. Future directions for research are focused on exploring global bifurcation behaviour, examining the behaviour of stability switches when changing several parameters, and furthering the scope of research to include systems composed of multiple interacting population groups. An additional area of interest for research is to examine the use of fractional or non-exponential memory kernels, which may allow for more complex dynamical patterns to develop and will provide further opportunity for exploration of the relationship between these types of memory structures and nonlocal evolution processes. In summary, these findings demonstrate how distributed state-dependent memory produces oscillatory behavior and offer a solid mathematical basis for more extensive future exploration of delayed driven dynamics.

\end{document}